\documentclass[10pt,a4paper]{article}
\usepackage[utf8]{inputenc}
\usepackage{amsmath}
\usepackage{amsfonts}
\usepackage{amssymb}
\usepackage{theorem}
\usepackage{mathrsfs}
\usepackage{cleveref}

\newtheorem{dfn}{Definition}[section]
\newtheorem{lem}[dfn]{Lemma}
\newtheorem{theor}[dfn]{Theorem}

\newtheorem{rem}[dfn]{Remark}
\newtheorem{cor}[dfn]{Corollary}
\newtheorem{ex}[dfn]{Example}

\newtheorem{pf}[dfn]{Proof.}
\crefname{section}{§}{§§}
\Crefname{section}{§}{§§}
\title{ $(C_\alpha,C_\beta)$-admissible functions in quasi-pseudometric type spaces. 
}

\author{ Ya\'e Ulrich Gaba, \\ gabayae2@gmail.com, \\ Department of Mathematics and Applied Mathematics,\\
University of Cape Town, Rondebosch 7701, South Africa.}

\begin{document}
\maketitle

\begin{abstract}
In this article, we give some fixed point results in left $K$-complete quasi-pseudometric type  spaces for self-mappings that are $(C_\alpha,C_\beta)$-admissible.

\end{abstract}
     
{\bf Keywords:} quasi-pseudometric; left(right)$K$-completeness; $(C_\alpha,C_\beta)$-admissible. 
     
\section{Introduction and Preliminaries}     
In the last few years the theory of quasi-metric spaces and other related structures such as quasi-normed cones and asymmetric normed linear spaces (see for instance \cite{ple6}) has been at the center of rigorous research activity, because such a theory provides a convenient framework in the study of several problems in theoretical computer science, approximation theory and convex analysis. The existence of fixed points and common fixed points of mappings satisfying certain contractive conditions in this setting, has also been discussed.

Recently, in \cite{eniola}, Eniola et al. discussed the newly introduced notion of \textit{quasi-pseudometric type  spaces} as a logical equivalent to metric type spaces- introduced by M. A. Khamsi\cite{khamsi}-when the initial distance-like function is not symmetric. 

The aim of this paper is to analyze the existence of fixed points for mapping defined on a left $K$-complete quasi-pseudometric type space $(X, D,K)$. The present results generalize another one already obtained by Gaba in \cite{rico}, where the author considered $(\alpha,\gamma)$-contractions.

\begin{dfn}
Let $X$ be a nonempty set, and let the function $D:X\times X \to [0,\infty)$ satisfy the following properties:
\begin{itemize}
\item[(D1)] $D(x,x)=0$ for any $x \in X$;
\item[(D2)] $D(x,y) \leq K \big( D(x,z_1)+D(z_1,z_2)+\cdots+D(x_n,y) \big)$ for any points $x,y,z_i\in X,\ i=1,2,\ldots, n$ and  some constant $K>0$.
\end{itemize}
Then $(X,D,K)$ is called a quasi-pseudometric type space.
Moreover, if $D(x,y)=0=D(y,x) \Longrightarrow x=y$, then $D$ is said to be a \textbf{$T_0$-quasi-pseudometric type space}. The latter condition is referred to as the $T_0$-condition.
\end{dfn}

\begin{rem} \hspace*{0.5cm}  
\begin{itemize}
\item Let $D$ be a quasi-pseudometric type on $X$, then the map $D^{-1}$
defined by $D^{-1}(x,y)=D(y,x)$ whenever $x,y \in X$ is also a quasi-pseudometric type on $X$, called the conjugate of $D$. We shall also denote $D^{-1}$ by $D^t$ or $\bar{D}$.
\item It is easy to verify that the function $D^s$ defined by $D^s:=D\vee D^{-1}$ , i.e. $D^s(x,y)=\max \{D(x,y),D(y,x)\}$  defines a \textbf{metric-type }(see \cite{khamsi}) on $X$ whenever $D$ is a $T_0$-quasi-pseudometric type.
\item If we substitute the property $(D1)$ by the following property 

$(D3)$ : $D(x,y)=0 \Longleftrightarrow x=y$, 

we obtain a $T_0$-quasi-pseudometric type space directly. 

\item  Moreover, for $K=1$, we recover the classical pseudometric, hence quasi-pseudmetric type spaces generalize quasi-pseudometrics. It is worth mentioning that if $(X,D,L)$ is a pseudometric type space, then for any $L \geq K$, $(X,D,L)$ is also a pseudometric type space. 
%Hence, in the sequel we shall denote $(X,D,\alpha)$ simply as $(X,D)$ when there is no confusion.
We give the following example to illustrate the above comment.

\begin{ex}
Let $X=\{ a,b,c \}$ and the mapping $D:X\times X \to [0,\infty)$ defined by $D(a,b)=D(c,b)=1/5,\ D(b,c)=D(b,a)=D(c,a)=1/4,\  D(a,c)=1/2$, $D(x,x)=0$ for any $x \in X$.
Since 
\[
\frac{1}{2} = D(a,c)> D(a,b)+D(b,c)=\frac{9}{20},
\]
then we conclude that $X$ is not a quasi-pseudometric space. Nevertheless, with $K=2$, it is very easy to check that $(X,D,2)$ is a quasi-pseudometric type space.
\end{ex}

\end{itemize}
\end{rem}

\begin{dfn} Let $(X, D, K)$ be a quasi-pseudometric type space.
The convergence of a sequence $(x_n)$ to $x$ with respect to $D$, called \textbf{$D$-convergence} or \textbf{left-convergence} and denoted by $x_n \overset{D}{\longrightarrow} x$,   
is defined in the following way
\begin{equation}
x_n \overset{D}{\longrightarrow} x \Longleftrightarrow D(x,x_n) \longrightarrow 0 .
\end{equation}

Similarly, the convergence of a sequence $(x_n)$ to $x$ with respect to $D^{-1}$, called \textbf{$D^{-1}$-convergence} or \textbf{right-convergence} and denoted by $x_n \overset{D^{-1}}{\longrightarrow} x$,

is defined in the following way
\begin{equation}
x_n \overset{D^{-1}}{\longrightarrow} x \Longleftrightarrow D(x_n,x) \longrightarrow 0 .
\end{equation}

Finally, in a quasi-pseudometric type space $(X, D, K)$, we shall say that a sequence $(x_n)$ \textbf{$D^s$-converges} to $x$ if it is both left and right convergent to $x$, and we denote it as $x_n \overset{D^{s}}{\longrightarrow} x$ or $x_n \longrightarrow x$ when there is no confusion.
Hence
\[
x_n \overset{D^{s}}{\longrightarrow} x \ \Longleftrightarrow \  x_n \overset{D}{\longrightarrow} x \ \text{ and }\ x_n \overset{D^{-1}}{\longrightarrow} x.
\]

\end{dfn}

\begin{dfn}
A sequence $(x_n)$ in a quasi-pseudometric type space $(X, D, K)$ is called
\begin{itemize}
\item[(a)] \textbf{left $K$-Cauchy } with respect to $D$ if for every $\epsilon >0$, there exists $n_0 \in \mathbb{N}$ such that 
$$ \forall \  n,k: n_0\leq k \leq n \quad D(x_k,x_n )< \epsilon ;$$
\item[(b)] \textbf{right $K$-Cauchy} with respect to $D$ if for every $\epsilon >0$, there exists $n_0 \in \mathbb{N}$ such that 
$$ \forall \  n,k: n_0\leq k \leq n \quad D(x_n,x_k )< \epsilon ;$$
\item[(c)] \textbf{ $D^s$-Cauchy} if for every $\epsilon >0$, there exists $n_0 \in \mathbb{N}$ such that 
$$ \forall n, k \geq n_0 \quad D(x_n,x_k )< \epsilon .$$
\end{itemize}
\end{dfn}

\begin{rem}\hspace*{0.5cm} 
\begin{itemize}
\item A sequence is left $K$-Cauchy with respect to $d$ if and only if it is right $K$-Cauchy with respect to $D^{-1}$. 
\item A sequence is $D^s$-Cauchy if and only if it is both left and right $K$-Cauchy.

\end{itemize}
\end{rem}

\begin{dfn}
 A quasi-pseudometric type space $(X, D, K)$ is called \textbf{left-complete} provided that any left $K$-Cauchy sequence is $D$-convergent.
\end{dfn}

\begin{dfn}
 A quasi-pseudometric type space $(X, D, K)$ is called \textbf{right-complete} provided that any right $K$-Cauchy sequence is $D$-convergent.
\end{dfn}

\begin{dfn}
 A $T_0$-quasi-pseudometric type space $(X, D, K)$ is called \textbf{bicomplete} provided that the metric type $D^s$ on $X$ is complete.
\end{dfn}

\section{Main Results}\label{sec:foo}

We begin by recalling the following.

\begin{dfn}(Compare \cite{rico})
Let $(X,D,K)$ be a quasi-pseudometric type space. A function $T:X\to X$ is called \textbf{$D$-sequentially continuous} or \textbf{left-sequentially continuous } if for any $D$-convergent sequence $(x_n)$ with $x_n \overset{D}{\longrightarrow} x $, the sequence $(Tx_n)$ $D$-converges to $Tx$, i.e. $Tx_n \overset{D}{\longrightarrow} Tx $.
\end{dfn}

We then introduce the following definition.
\begin{dfn}
Let $(X,D,K)$ be a quasi-pseudometric type space, $f:X\to X$ and $\alpha,\beta: X \times X\to [0,\infty)$ be mappings and $C_\alpha>0, C_\beta\geq 0$. We say that $T$ is $(C_\alpha,C_\beta)$-admissible with respect to $K$ if the following conditions hold: 
\begin{enumerate}
\item[(C1)] $\alpha(x,y)\geq C_\alpha \Longrightarrow \alpha(fx,fy)\geq C_\alpha,$ whenever $x,y\in X$;
\item[(C2)]$\beta(x,y)\leq C_\beta \Longrightarrow \beta(fx,fy)\leq C_\beta,$ whenever $x,y\in X$;
\item[(C3)] $0\leq C_\beta/C_\alpha < 1/K$.
\end{enumerate}
\end{dfn}

We begin by the following lemma.
\begin{lem}\label{lemme1}
Let $(X,D,K)$ be a quasi-pseudometric type space and let $(x_n)$ be a sequence in $X$. Then
\[
D(x_0,x_n) \leq K D(x_0,x_1)+\cdots + K^{n-1}D(x_{n-2},x_{n-1})+K^{n-1}D(x_{n-1},x_n).
\]
\end{lem}
From Lemma \ref{lemme1}, we deduce the following lemma.

\begin{lem}\label{cochi}(Compare \cite[Lemma 38]{eniola})
Let $(X,D,K)$ be a quasi-pseudometric type space and let $(x_n)$ be a sequence in $X$ such that

\[
D(x_n,x_{n+1}) \leq \lambda D(x_{n-1},x_n)\qquad \text{ for all } \ n\geq 0,
\]
for some $0<\lambda<1/K$. Then $(x_n)$ is a left $K$-Cauchy sequence.
\end{lem}

Similarly, 

\begin{lem}
Let $(X,D,\alpha)$ be a quasi-pseudometric type space and let $(x_n)$ be a sequence in $X$ such that

\[
D(x_{n+1},x_n) \leq \lambda D(x_n,x_{n-1}) \qquad \text{ for all } \ n\geq 0,
\]
for some $0<\lambda<1/K$. Then $(x_n)$ is a right $K$-Cauchy sequence.
\end{lem}

We now state our main fixed point theorem.

\begin{theor}\label{fix1}
Let $(X,D,K)$ be a Hausdorff left $K$-complete $T_0$-quasi-pseudometric type space. Suppose that $f:X\to X$ is $(C_\alpha,C_\beta)$-admissible with respect to $K$. Assume that
\begin{equation}\label{cond}
\alpha(x,y)D(fx,fy) \leq \beta(x,y)D(x,y) \qquad \text{for all } \ \ x,y \in X.
\end{equation}

 If the following conditions hold:

\begin{enumerate}
\item[(i)] $f$ is $D$-sequentially continuous;
\item[(ii)] there exists $x_0\in X$ such that $\alpha(x_0,fx_0)\geq C_\alpha$ and
$\beta(x_0,fx_0)\leq C_\beta$.
\end{enumerate}
Then $f$ has a fixed point.
\end{theor}

\begin{pf}
Let $x_0\in X$ such that $\alpha(x_0,fx_0)\geq C_\alpha$ and $\beta(x_0,fx_0)\leq C_\beta$. Define the sequence $(x_n)$ by $x_n=f^nx_0=fx_{n-1}$. Without loss of generality, we can always assume that $x_n\neq x_{n+1}$ for all $n\in\mathbb{N}$, since if $x_{n_0}=x_{n_0+1}$ for some $n_0\in \mathbb{N},$ the proof is complete. Since $f$ is  $(C_\alpha,C_\beta)$-admissible with respect to $K$ and $\alpha(x_0,fx_0)= (x_0,x_1)\geq C_\alpha,$
we deduce that $\alpha(x_1,x_2)= (fx_0,fx_1)\geq C_\alpha.$ By continuing this process, we get that
$\alpha(x_n,x_{n+1})\geq C_\alpha,$ for all $n\geq 0$. Similarly, we establish that $\beta(x_n,x_{n+1})\leq C_\beta,$ for all $n\geq 0$. Using \eqref{cond}, we get 

\begin{align*}
  C_\alpha D(x_n,x_{n+1}) & \leq  \alpha(x_{n-1},x_{n}) D(x_n,x_{n+1}) \\
                          & \leq  \beta(x_{n-1},x_{n})  D(x_{n-1},x_{n}) \\
                          & \leq C_\beta D(x_{n-1},x_n),
\end{align*}
and hence 
\[
D(x_n,x_{n+1}) \leq \frac{C_\beta}{C_\alpha} D(x_{n-1},x_n) \qquad \text{ for all } \ n\geq 1.
\]

By Lemma \ref{cochi}, since $0\leq  \frac{C_\beta}{C_\alpha} <1/K,$ we derive that $(x_n)$ is a left $K$-Cauchy sequence. Since $(X,d)$ is left $K$-complete and $T$ $D$-sequentially continuous, there exists $x^*$ such that $x_n \overset{D}{\longrightarrow} x^*$ and $x_{n+1}\overset{D}{\longrightarrow} fx^*$. Since $X$ is Hausdorff, $x^*=fx^*$.
\end{pf}

\begin{cor}
Let $(X,D,K)$ be a Hausdorff left $K$-complete $T_0$-quasi-pseudometric type space. Suppose that $f:X\to X$ is $(C_\alpha,C_\beta)$-admissible with respect to $K$. Assume that
\begin{equation}
\alpha(x,y)D(fx,fy) \leq \beta(x,y)D(x,y) \qquad \text{for all } \ \ x,y \in X.
\end{equation}

 If the following conditions hold:

\begin{enumerate}
\item[(i)] $f$ is $D^{-1}$-sequentially continuous;
\item[(ii)] there exists $x_0\in X$ such that $\alpha(fx_0,x_0)\geq C_\alpha$ and 
$\beta(fx_0,x_0)\leq C_\beta$.
\end{enumerate}
Then $f$ has a fixed point.
\end{cor}

\begin{cor}
Let $(X,D,K)$ be a bicomplete $T_0$-quasi-pseudometric type space. Suppose that $f:X\to X$ is $(C_\alpha,C_\beta)$-admissible with respect to $K$. Assume that
\begin{equation}
\alpha(x,y)D(fx,fy) \leq \beta(x,y)D(x,y) \qquad \text{for all } \ \ x,y \in X.
\end{equation}
If the following conditions hold:

\begin{enumerate}
\item[(i)] $f$ is $D^s$-sequentially continuous;
\item[(ii)] there exists $x_0\in X$ such that $\alpha(x_0,fx_0)\geq C_\alpha$ and
$\beta(x_0,fx_0)\leq C_\beta$;
\item[(iii)] the functions $\alpha$ and $\beta$ are symmetric,  i.e. $\alpha(a,b)=\alpha(b,a)$ for any $a,b \in X$ and $\beta(a,b)=\beta(b,a)$ for any $a,b \in X$.
\end{enumerate}
Then $f$ has a fixed point.
\end{cor}

\begin{pf}
Following the proof of Theorem \ref{fix1}, it is clear that the sequence $(x_n)$ in $X$ defined by  $x_{n+1}=fx_n$ for all $n=0,1,2,\cdots$ is  $D^s$-Cauchy. Since $(X,D^s)$ is complete and $f$ $D^s$-sequentially continuous, there exists $x^*$ such that $x_n \overset{D^s}{\longrightarrow} x^*$ and $x_{n+1}\overset{D^s}{\longrightarrow} fx^*$. Since $(X,D^s)$ is Hausdorff, $x^*=fx^*$.
\end{pf}

\begin{rem}\label{remark}
In fact, we do not need $\alpha$ and $\beta$ to be symmetric. It is enough, for the result to be true, to have a  point $x_0 \in X$ for which $\alpha(x_0,fx_0)\geq C_\alpha$, $\alpha(fx_0,x_0)\geq C_\alpha$, $\beta(x_0,fx_0)\leq C_\beta$ and $\beta(fx_0,x_0)\leq C_\beta$. Otherwise stated, it is enough to have a  point $x_0 \in X$ for which $\min\{\alpha(x_0,fx_0),\alpha(fx_0,x_0)\} \geq C_\alpha$ and $\min \{\beta(x_0,fx_0),\beta(fx_0,x_0\}\leq C_\beta$.
\end{rem}

We give the following results which are in fact consequences of the Theorem \ref{fix1}.

\begin{theor}
Let $(X,D,K)$ be a Hausdorff left $K$-complete $T_0$-quasi-pseudometric type space. Suppose that $f:X\to X$ is $(C_\alpha,C_\beta)$-admissible with respect to $K$. Assume that
\begin{equation}\label{conds}
\alpha(x,y)D(fx,fy) \leq \beta(x,y)D(x,y) \qquad \text{for all } \ \ x,y \in X.
\end{equation}

 If the following conditions hold:

\begin{enumerate}
\item[(i)] there exists $x_0\in X$ such that $\alpha(x_0,fx_0)\geq C_\alpha$ and
$\beta(x_0,fx_0)\leq C_\beta$;
\item[(ii)] if $(x_n)$ is a sequence in $X$ such that $\alpha(x_n,x_{n+1})\geq C_\alpha$  and $\beta(x_n,x_{n+1})\leq C_\beta$ for all $n=1,2,\cdots$ and $x_n \overset{D}{\longrightarrow} x$, then there exists a subsequence $(x_{n(k)})$ of $(x_n)$ such that $\alpha(x,x_{n(k)})\geq C_\alpha$ and $\beta(x,x_{n(k)})\leq C_\beta$ for all $k$.
\end{enumerate}
Then $f$ has a fixed point.
\end{theor}

\begin{pf}
Following the proof of Theorem \ref{fix1}, we know that the sequence $(x_n)$ defined by 
$x_{n+1}=fx_n$ for all $n=0,1,2,\cdots$ $D$-converges to some $x^*$ and satisfies $\alpha(x_n,x_{n+1})\geq C_\alpha$ and $\beta(x_n,x_{n+1})\leq C_\beta $ for $n\geq 1$. From the condition $(ii)$, we know there exists a subsequence $(x_{n(k)})$ of $(x_n)$ such that $\alpha(x^*,x_{n(k)})\geq C_\alpha$ and $\beta(x^*x_{n(k)})\leq C_\beta$ for all $k$. Since $f$ is satisfies \eqref{conds}, we get 

\begin{align*}
  C_\alpha D(fx^*,x_{n(k)+1}) = C_\alpha D(fx^*,fx_{n(k)})& \leq  \alpha(x^*,x_{n(k)}) D(fx^*,fx_{n(k)}) \\
                          & \leq   \beta(x^*,x_{n(k)})  D(x^*,x_{n(k)}) \\
                          & \leq C_\beta  D(x^*,x_{n(k)})
\end{align*}

Letting $k\to \infty$, we obtain $D(fx^*,x_{n(k)+1}) \to 0$. Since $X$ is Hausdorff, we have that $fx^*=x^*$.
This completes the proof.
\end{pf}

\begin{cor}
Let $(X,D,K)$ be a Hausdorff right $K$-complete $T_0$-quasi-pseudometric type space. Suppose that $f:X\to X$ is $(C_\alpha,C_\beta)$-admissible with respect to $K$. Assume that
\begin{equation}
\alpha(x,y)D(fx,fy) \leq \beta(x,y)D(x,y) \qquad \text{for all } \ \ x,y \in X.
\end{equation}

 If the following conditions hold:

\begin{enumerate}
\item[(i)] there exists $x_0\in X$ such that $\alpha(fx_0,x_0)\geq C_\alpha$ and
$\beta(fx_0,x_0)\leq C_\beta$;
\item[(ii)] if $(x_n)$ is a sequence in $X$ such that $\alpha(x_{n+1},x_n)\geq C_\alpha$  and $\beta(x_{n+1},x_n)\leq C_\beta$ for all $n=1,2,\cdots$ and $x_n \overset{D^{-1}}{\longrightarrow} x$, then there exists a subsequence $(x_{n(k)})$ of $(x_n)$ such that $\alpha(x_{n(k)},x)\geq C_\alpha$ and $\beta(x_{n(k)},x)\leq C_\beta$ for all $k$.
\end{enumerate}
Then $f$ has a fixed point.
\end{cor}

\begin{cor}
Let $(X,D,K)$ be a bicomplete $T_0$-quasi-pseudometric type space. Suppose that $f:X\to X$ is $(C_\alpha,C_\beta)$-admissible with respect to $K$. Assume that
\begin{equation}
\alpha(x,y)D(fx,fy) \leq \beta(x,y)D(x,y) \qquad \text{for all } \ \ x,y \in X.
\end{equation}

 If the following conditions hold:
 \begin{enumerate}
\item[(i)] there exists $x_0\in X$ such that $\alpha(x_0,fx_0)\geq C_\alpha$ and $\beta(x_0,fx_0)\leq C_\beta$;

\item[(ii)] if $(x_n)$ is a sequence in $X$ such that $\alpha(x_m,x_n)\geq C_\alpha$  and $\beta(x_m,x_n)\leq C_\beta$ for all $n,m \in \mathbb{N}$ and $x_n \overset{D^s}{\longrightarrow} x$, then there exists a subsequence $(x_{n(k)})$ of $(x_n)$ such that $\alpha(x,x_{n(k)})\geq C_\alpha$ and $\beta(x,x_{n(k)})\leq C_\beta$ for all $k$;
\item[(iii)] the functions $\alpha,\beta$ are symmetric.
\end{enumerate}
Then $f$ has a fixed point.
\end{cor}

In regard of Remark \eqref{remark}, another variant of the above corollary can be stated as follows:
\begin{cor}
Let $(X,D,K)$ be a bicomplete $T_0$-quasi-pseudometric type space. Suppose that $f:X\to X$ is $(C_\alpha,C_\beta)$-admissible with respect to $K$. Assume that
\begin{equation}
\alpha(x,y)D(fx,fy) \leq \beta(x,y)D(x,y) \qquad \text{for all } \ \ x,y \in X.
\end{equation}

If the following conditions hold:
\begin{enumerate}
\item[(i)] there exists $x_0\in X$ such that $\min \{\alpha(fx_0,x_0), \alpha(x_0,fx_0)\}\geq C_\alpha$ and $\min \{\beta(fx_0,x_0),\beta(x_0,fx_0)\}\leq C_\beta$;

\item[(ii)] if $(x_n)$ is a sequence in $X$ such that $\min\{\alpha(x_m,x_n),\alpha(x_n,x_m)\}\geq C_\alpha$  and $\min\{\beta(x_m,x_n),\beta(x_n,x_m)\}\leq C_\beta$ for all $n,m \in \mathbb{N}$ and $x_n \overset{D^s}{\longrightarrow} x$, then there exists a subsequence $(x_{n(k)})$ of $(x_n)$ such that $\min\{\alpha(x,x_{n(k)}),\alpha(x_{n(k)},x) \}\geq C_\alpha$ and $\min\{\beta(x,x_{n(k)}),\beta(x_{n(k)},x) \}\leq C_\beta$ for all $k$.
\end{enumerate}
Then $f$ has a fixed point.
\end{cor}

We conclude this paper by the following theorem.

\begin{theor}\label{fix2}
Let $(X,D,K)$ be a bicomplete $T_0$-quasi-pseudometric space. Suppose that $f:X\to X$ is $(C_\alpha,C_\beta)$-admissible with respect to $K$. Assume that
\begin{equation}\label{condi}
\alpha(x,y)D^s(fx,fy) \leq \beta(x,y)D^s(x,y) \qquad \text{for all } \ \ x,y \in X.
\end{equation}

 If the following conditions hold:

\begin{enumerate}
\item[(i)] there exists $x_0\in X$ such that $\alpha(x_0,fx_0)\geq C_\alpha$ and
$\beta(x_0,fx_0)\leq C_\beta$;
\item[(ii)] if $(x_n)$ is a sequence in $X$ such that $\alpha(x_{n+1},x_n)\geq C_\alpha$  and $\beta(x_{n+1},x_n)\leq C_\beta$ for all $n=0,1,2,\cdots$ and $x_n \overset{D^{s}}{\longrightarrow} x$, then $\alpha(x_n,x)\geq C_\alpha$ and $\beta(x_n,x)\leq C_\beta$ for all $n\in \mathbb{N}$.
\end{enumerate}
Then $f$ has a fixed point.
\end{theor}

\begin{pf}
Let $x_0\in X$ such that $\alpha(x_0,fx_0)\geq C_\alpha$ and ;
$\beta(x_0,fx_0)\leq C_\beta$. Define a sequence $(x_n)$ by $x_n=f^nx_0=fx_{n-1}$. Without loss of generality, we can always assume that $x_n\neq x_{n+1}$ for all $n\in\mathbb{N}$, since if $x_{n_0}=x_{n_0+1}$ for some $n_0\in \mathbb{N},$ the proof is complete. Following the proof of Theorem \ref{fix1}, we see that $(x_n)$ is a $D^s$-Cauchy sequence $\alpha(x_{n+1},x_n)\geq C_\alpha$  and $\beta(x_{n+1},x_n)\leq C_\beta$ for all $n=0,1,2,\cdots$. Since $(X,D^s)$ is complete, there exists $x^*$ such that $x_n \overset{D^s}{\longrightarrow} x^*$ and. Now, using the contractive condition \eqref{condi} and condition $(ii)$, we deduce that

\begin{align*}
   D^s(x^*,fx^*) & \leq K\left[  D^s(x^*,x_{n+1})+ \frac{C_\alpha}{C_\alpha} D^s(x_{n+1},fx^*))\right]\\
               & \leq K D^s(x^*,x_{n+1})  + \frac{K}{C_\alpha}\alpha(x_n,x^*) D^s(fx_{n},fx^*)\\
               & \leq K D^s(x^*,x_{n+1})  + \frac{K}{C_\alpha}\beta(x_n,x^*) D^s(x_{n},x^*))\\
               & \leq K D^s(x^*,x_{n+1})  + \frac{K\beta}{C_\alpha} D^s(x_{n},x^*).
\end{align*}
 Letting $n\to \infty$, we obtain that $D^s(x^*,fx^*)\leq 0,$ that is $x^*=fx^*.$ This complete the proof.

\end{pf}

\end{document}